\documentclass[12pt,reqno,a4paper]{amsart}
\usepackage[a4paper,margin=1in]{geometry}

\title[Ambient Chow rings of Debarre--Voisin fourfolds]{Ambient Chow Rings of Debarre--Voisin Fourfolds and the Beauville--Voisin Conjecture}
\author{Ryo Akazawa}

\address{
Department of Mathematics, Faculty of Science and Engineering,
Waseda University,
3-4-1 Ohkubo, Shinjuku, Tokyo 169-8555, Japan}
\email{rakazawa11235@ruri.waseda.jp}

\usepackage{amsmath,amssymb}
\usepackage{mathtools}
\usepackage{newtxtext}
\usepackage[nosymbolsc]{newtxmath}
\usepackage[mathscr]{eucal}

\usepackage{graphicx}
\usepackage{array,booktabs}
\usepackage{listings}

\usepackage{xcolor}
\usepackage{hyperref}
\hypersetup{%
 colorlinks=true,%
 linkcolor=blue!75!black,%
 citecolor=red!70!black,%
 urlcolor=teal!85!black%
}

\theoremstyle{plain}
\newtheorem{thm}{Theorem}[section]
\newtheorem{lem}[thm]{Lemma}
\newtheorem{prop}[thm]{Proposition}
\newtheorem{cor}[thm]{Corollary}
\newtheorem{conj}[thm]{Conjecture}

\newtheorem*{thm*}{Theorem}
\newtheorem*{cor*}{Corollary}

\theoremstyle{definition}
\newtheorem{dfn}[thm]{Definition}

\theoremstyle{remark}

\newtheorem{rem}[thm]{Remark}

\DeclareMathOperator{\Pic}{Pic}

\DeclareMathOperator{\Hom}{Hom}

\DeclareMathOperator{\ch}{ch}

\DeclareMathOperator{\Gr}{Gr}

\DeclareMathOperator{\Ker}{Ker}

\newcommand{\C}{\mathbb{C}}

\newcommand{\Q}{\mathbb{Q}}

\newcommand{\cE}{\mathcal E}

\makeatletter
\def\@secnumfont@section{}
\def\@secnumfont@subsection{\protect\bfseries}
\def\@seccntformat#1{%
  \protect\textup{\protect\@secnumfont
    \csname @secnumfont@#1\endcsname
    \csname the#1\endcsname
    \protect\@secnumpunct
  }%
}
\makeatother

\numberwithin{equation}{section}

\providecommand{\Q}{\mathbb Q}
\providecommand{\C}{\mathbb C}

\providecommand{\CH}{\mathrm{CH}}
\providecommand{\Gr}{\operatorname{Gr}}

\providecommand{\cU}{\mathcal U}
\providecommand{\cQ}{\mathcal Q}
\providecommand{\cK}{\mathcal K}
\providecommand{\cE}{\mathcal E}
\providecommand{\cO}{\mathcal O}
\providecommand{\cl}{\mathrm{cl}}

\providecommand{\Imop}{\operatorname{Im}}
\providecommand{\Ker}{\operatorname{Ker}}

\providecommand{\Pic}{\operatorname{Pic}}

\providecommand{\ch}{\operatorname{ch}}
\providecommand{\Hom}{\operatorname{Hom}}

\providecommand{\colonequals}{\coloneqq}

\begin{document}

\begin{abstract}
We prove the Beauville--Voisin conjecture for a very general
Debarre--Voisin hyperk\"ahler fourfold \(Y\subset\Gr(6,10)\)
associated with a trivector \(\sigma\) on a ten-dimensional complex
vector space \(V_{10}\). Namely, the cycle class map is injective on
the subalgebra of the Chow ring generated by divisor classes and the
Chern classes of the tangent bundle.

\end{abstract}

\maketitle
\setcounter{tocdepth}{1}
\tableofcontents

\section{Introduction}
\label{sec:introduction}

The \emph{Chow ring} of a smooth complex projective variety \(X\)
consists of algebraic cycles modulo rational equivalence, equipped with
the intersection product. It is related to the cohomology ring by the
\emph{cycle class map}
\[
\cl:\CH^i(X)\longrightarrow H^{2i}(X,\Q).
\]
Throughout this paper, Chow groups are taken with rational coefficients,
and \(H^*(X,\Q)\) denotes singular cohomology.

The Chow ring contains rich geometric information, but is often difficult
to describe. Mumford's work on zero-cycles on surfaces with holomorphic
two-forms \cite{Mumford} illustrates how large the degree-zero subgroup
\(\Ker(\deg:\CH_0(X)\to\Q)\) can be. For \emph{K3 surfaces}, however,
geometrically distinguished classes satisfy strong constraints, as shown
by the following theorem of Beauville and Voisin.

\begin{theorem}[Beauville--Voisin, {\cite[Theorem 1]{BeauvilleVoisin}}]\label{thm:bv-k3}
Let \(S\) be a projective K3 surface. There is a class
\(o_S\in\CH^2(S)\) of degree one such that, for any divisor classes
\(D,D'\),
\[
DD'=\deg(DD')o_S,\qquad c_2(T_S)=24o_S.
\]
The class \(o_S\) is represented by any point on any rational curve on \(S\).
\end{theorem}

Consequently, the subring generated by divisor classes and \(c_2(T_S)\) is
\(\Q[S]\oplus\CH^1(S)\oplus\Q o_S\), and its cycle class map is injective.

A projective \emph{irreducible holomorphic symplectic variety} (IHS),
also called a \emph{hyperk\"ahler variety} (HK), is a simply connected
smooth complex projective variety whose space \(H^0(X,\Omega_X^2)\)
is generated by a nondegenerate holomorphic two-form.
Such varieties are natural higher-dimensional analogues of K3 surfaces.
Beauville conjectured that the cycle class map is injective on the
subring generated by divisor classes \cite{BeauvilleSplitting}.
Voisin extended this conjecture
to include the Chern classes of the tangent bundle
\cite[Conjecture 1.3]{VoisinChow}.

\begin{conjecture}[Beauville--Voisin]\label{conj:bv}
Let \(X\) be a projective irreducible holomorphic symplectic variety, and let
\[
R^*(X)\colonequals \langle\CH^1(X),c_i(T_X)\rangle_{\Q}
\subset\CH^*(X)
\]
be the \(\Q\)-subalgebra generated by divisor classes and the Chern
classes of the tangent bundle. Then \(\cl|_{R^*(X)}\) is injective.
\end{conjecture}

Voisin proved this conjecture for Hilbert schemes of points on K3
surfaces under certain hypotheses, and for the Fano variety of lines
on a smooth cubic fourfold \cite{VoisinChow}.
Yin proved it for Hilbert schemes of arbitrary length on K3 surfaces
with finite-dimensional Chow motive \cite{YinFinite}.
Fu established it for generalized Kummer varieties \cite{FuKummer}.

Further known cases include smooth double EPW sextics
\cite{LaterveerEPW} and very general
Lehn--Lehn--Sorger--van Straten eightfolds
\cite{FLVFranchettaII}.
More recently, Mazzanti proved the conjecture for double EPW
quartics in a preprint \cite{MazzantiEPWQuartics}.

We study \emph{Debarre--Voisin fourfolds}. Let \(V=V_{10}\) be a
ten-dimensional complex vector space, and let
\(0\ne\sigma\in\bigwedge^3V^\vee\). Set
\[
\begin{split}
Y_\sigma&\colonequals \{W\in\Gr(6,V)\mid\sigma|_W=0\},\\
F_\sigma&\colonequals \{K\in\Gr(3,V)\mid\sigma|_K=0\}.
\end{split}
\]
Here \(G_k\colonequals\Gr(k,V)\) is the Grassmannian of
\(k\)-dimensional subspaces of \(V\). If \(Y=Y_\sigma\) is a smooth
four-dimensional zero scheme, then \(F_\sigma\) is also smooth, and
\(Y\) is a hyperk\"ahler fourfold deformation equivalent to the Hilbert
square of a K3 surface \cite{DebarreVoisin,DHOV}.
The incidence correspondence \(G_\sigma\) between \(Y\) and the
twenty-dimensional Pl\"ucker hyperplane section \(F_\sigma\) relates
their Hodge structures. Voisin proved a vanishing theorem for the Chow
groups of \(F_\sigma\) \cite{VoisinGeneralized}, and Laterveer studied
its Chow ring and generically defined cycles \cite{LaterveerDV}.

Let \(i:Y\hookrightarrow G_6\) be the inclusion, and let \(\cQ\) be the
universal quotient bundle on \(G_6\).
For its restriction \(\cQ_Y=i^*\cQ\), set
\[
q_j\colonequals c_j(\cQ_Y),\qquad h\colonequals q_1,\qquad z\colonequals 44q_2-19h^2.
\]
The class \(h\) is the first Chern class of the Pl\"ucker polarization.
We call the subring of classes restricted from the Grassmannian the
\emph{ambient Chow ring}:
\[
\CH^*_{\mathrm{amb}}(Y)\colonequals \Imop\left(i^*:\CH^*(G_6)\to\CH^*(Y)\right).
\]
Our main result is the following.

\begin{theorem}\label{thm:main}
Suppose that \(Y=Y_\sigma\) is a smooth four-dimensional zero scheme,
that is, a Debarre--Voisin fourfold.
Then the cycle class map is injective on \(\CH^*_{\mathrm{amb}}(Y)\),
and there is an isomorphism of graded rings
\[
\Q[H,\zeta]/(H\zeta,\zeta^2-11H^4)
\xrightarrow{\ \sim\ }\CH^*_{\mathrm{amb}}(Y),
\qquad H\longmapsto h,\quad\zeta\longmapsto z,
\]
where \(\deg H=1\) and \(\deg\zeta=2\).
In particular, \(\CH^3_{\mathrm{amb}}(Y)=\Q h^3\) and
\(\CH^4_{\mathrm{amb}}(Y)=\Q h^4\).
\end{theorem}

\begin{corollary}\label{cor:bv}
Every smooth Debarre--Voisin fourfold of Picard number one satisfies
Conjecture~\ref{conj:bv}. In particular, a very general \(Y_\sigma\)
satisfies this conjecture.
\end{corollary}

Here \emph{very general} means outside a countable union of proper
closed subsets of the parameter space.

\subsection*{Organization of the paper}

Section~\ref{sec:preliminaries} introduces the bundles and push-forward
formulas used in the proof. In Section~\ref{sec:relations}, we first
deduce \(3h^3-8hq_2+5q_3=0\) from the vanishing of \(c_3(T_Y)\).
We then use the incidence correspondence between \(Y\) and \(F_\sigma\)
to obtain two further relations from the vanishing theorem for the
Chow groups of \(F_\sigma\). The remaining relation comes from a
section of a rank-twelve bundle on this incidence variety, constructed
from \(\sigma\). Smoothness of \(Y\) implies that this section is
nowhere vanishing, so the top Chern class of the bundle is zero.
We compute its push-forward using the Gysin formula for Grassmann
bundles.
In Section~\ref{sec:ambient-ring}, we use these relations to give a
presentation of the ambient ring, and prove the main theorem by
checking linear independence in cohomology using intersection products.

In Section~\ref{sec:applications}, we compute \(c_2(T_Y)\) and
\(c_4(T_Y)\), and deduce Corollary~\ref{cor:bv} from the resulting
Chow relations and the main theorem. Finally, we apply
Theorem~\ref{thm:main} to the universal family parametrized by
\([\sigma]\). Classes restricted from the total space of this family
to a fiber are precisely the ambient classes, and hence vanish in the
Chow group whenever their cohomology classes vanish. This is the
Franchetta property for this family, in the framework of the
generalized Franchetta conjecture for hyperk\"ahler families
\cite{Franchetta}.

\subsection*{Acknowledgments}

The author is sincerely grateful to his advisor, Professor Yasunari
Nagai, for helpful advice and valuable discussions, and for his
encouragement.

\section{Preliminaries}
\label{sec:preliminaries}

Throughout the rest of the paper, \(Y=Y_\sigma\) is assumed to be a
smooth Debarre--Voisin fourfold. All Chow groups have rational
coefficients. For a smooth variety \(X\), we write
\(\CH_k(X)=\CH^{\dim X-k}(X)\) for the Chow group of
\(k\)-dimensional cycles.

\subsection{Notation and tautological bundles}

Let \(G_6=\Gr(6,V)\), and let \(i:Y\hookrightarrow G_6\) be the
inclusion. The tautological sequence on \(G_6\) is
\[
0\longrightarrow\cU\longrightarrow V\otimes\cO_{G_6}
\longrightarrow\cQ\longrightarrow0,
\]
where \(\cU\) and \(\cQ\) have ranks six and four, respectively.
For their restrictions, set
\[
\cU_Y\colonequals i^*\cU,\qquad \cQ_Y\colonequals i^*\cQ,\qquad
u_j\colonequals c_j(\cU_Y^\vee),\qquad q_j\colonequals c_j(\cQ_Y),\qquad h\colonequals q_1.
\]
The Whitney formula gives
\begin{equation}\label{eq:whitney-y}
c(\cU_Y)c(\cQ_Y)=1.
\end{equation}
Since \(c_j(\cU_Y)=(-1)^j u_j\), comparison of the components of
codimensions one through four in \eqref{eq:whitney-y} gives
\begin{equation}\label{eq:chern-components}
\begin{aligned}
-u_1+h&=0,\\
u_2-hu_1+q_2&=0,\\
-u_3+hu_2-q_2u_1+q_3&=0,\\
u_4-hu_3+q_2u_2-q_3u_1+q_4&=0.
\end{aligned}
\end{equation}
Solving \eqref{eq:chern-components} successively yields
\begin{equation}\label{eq:u-q}
\begin{aligned}
u_1&=h,& u_2&=h^2-q_2,\\
u_3&=h^3-2hq_2+q_3,&
u_4&=h^4-3h^2q_2+2hq_3+q_2^2-q_4.
\end{aligned}
\end{equation}

We use the conventions \(q_0=1\) and \(q_j=0\) for \(j<0\) or
\(j>4\). The Chow ring \(\CH^*(G_6)\) is generated by the classes
\(c_j(\cQ)\). Thus the ambient groups
\(\CH^*_{\mathrm{amb}}(Y)=\Imop(i^*)\) are spanned in each
codimension by the following classes:
\begin{equation}\label{eq:generators}
\begin{array}{c|l}
\multicolumn{1}{c}{\text{codim}}&\\
\hline
0&1\\
1&h\\
2&h^2,q_2\\
3&h^3,hq_2,q_3\\
4&h^4,h^2q_2,hq_3,q_2^2,q_4.
\end{array}
\end{equation}
The class \(h\) is the first Chern class of the ample Pl\"ucker
polarization, so \(d\colonequals \int_Yh^4>0\). It follows that the cohomology
class of \(h^j\) is nonzero for \(0\le j\le4\): otherwise,
multiplication by \(h^{4-j}\) would give \(\int_Yh^4=0\).

\subsection{The defining section and the incidence correspondence}

Let \(E\colonequals \bigwedge^3\cU^\vee\), and let \(s_\sigma\) be the section
of \(E\) defined by \(\sigma\). We have
\(\operatorname{rank}E=\binom63=20\) and
\(\dim G_6=6(10-6)=24\).

At a point \(W\in Y\), the tangent space to the Grassmannian is
\(T_{G_6,W}=\Hom(W,V/W)\). Write
\(d\sigma_W\colonequals (ds_\sigma)_W:\Hom(W,V/W)\to\bigwedge^3W^\vee\)
for the differential of the section. By the definition of the zero
scheme, \(T_{Y,W}=\Ker(d\sigma_W)\). Since \(Y\) is smooth of
dimension four, the differential has rank \(24-4=20\). Its target
also has dimension twenty, so \(d\sigma_W\) is surjective.

In a local trivialization of \(E\), write
\(s_\sigma=(f_1,\ldots,f_{20})\). The differentials
\(df_1,\ldots,df_{20}\) are linearly independent at each \(W\in Y\).
Hence \(f_1,\ldots,f_{20}\) form part of a regular system of
parameters in the local ring of \(G_6\), and in particular form a
regular sequence. Thus \(s_\sigma\) is a regular section and
\(N_{Y/G_6}\simeq E|_Y\). Using
\(T_{G_6}\simeq\cU^\vee\otimes\cQ\), we obtain the normal sequence
\begin{equation}\label{eq:normal}
0\longrightarrow T_Y\longrightarrow
\cU_Y^\vee\otimes\cQ_Y\xrightarrow{\,d\sigma\,}
\bigwedge^3\cU_Y^\vee\longrightarrow0.
\end{equation}

For \(f\in\Hom(W,V/W)\), choose a lift \(\widetilde f:W\to V\).
For \(w_1,w_2,w_3\in W\), the differential is given by
\begin{align}
d\sigma_W(f)(w_1,w_2,w_3)={}&\sigma(\widetilde f(w_1),w_2,w_3)
+\sigma(w_1,\widetilde f(w_2),w_3)\notag\\
&+\sigma(w_1,w_2,\widetilde f(w_3)).\label{eq:differential}
\end{align}
Indeed, after choosing a complement to \(W\), restrict \(\sigma\)
to the graph of \(tf\) and take the coefficient of \(t\). This gives
\eqref{eq:differential}. The expression is independent of the lift,
since \(\sigma|_W=0\).

Let \(G_3=\Gr(3,V)\), with \(F_\sigma\subset G_3\).

\begin{lemma}[{\cite[Proposition 4.4(a)]{DHOV}}]\label{lem:smooth-fsigma}
The variety \(F_\sigma\) is a smooth Pl\"ucker hyperplane section.
\end{lemma}
\begin{proof}
By the cited proposition, a singular point \(K\) of \(F_\sigma\)
would be contained in some \(W\in Y\) with \(\dim T_{Y,W}>4\).
This contradicts the assumption that \(Y\) is smooth of dimension four.
\end{proof}

Let \(j:F_\sigma\hookrightarrow G_3\) be the inclusion.
The incidence correspondence of Debarre and Voisin
\cite[Section 2]{DebarreVoisin} is
\[
G_\sigma\colonequals \{(K,W)\in F_\sigma\times Y\mid K\subset W\}
\simeq\Gr_Y(3,\cU_Y).
\]
Indeed, \(K\subset W\) and \(\sigma|_W=0\) imply \(K\in F_\sigma\),
and the relative Grassmannian parametrizes exactly these pairs of
subspaces. We denote the projections by
\[
Y_\sigma\xleftarrow{\ q\ }G_\sigma\xrightarrow{\ p\ }F_\sigma,
\qquad q(K,W)=W,\quad p(K,W)=K.
\]
The relative dimension of \(q\) is \(\dim\Gr(3,6)=9\).
Since \(Y\), \(F_\sigma\), and \(G_\sigma\) are smooth, \(p^*\)
is defined as the intersection-theoretic pullback.

\subsection{Schubert classes and push-forward formulas}

A \emph{partition} \(\lambda=(\lambda_1\ge\lambda_2\ge\cdots\ge0)\)
is a sequence of nonnegative integers with only finitely many nonzero
terms. Its size is \(|\lambda|=\sum\lambda_j\), and its length
\(\ell(\lambda)\) is the number of positive terms.
A \emph{Schubert variety} is a subvariety defined by intersection conditions
with a fixed complete flag; its Chow class is called a \emph{Schubert class}.
We write \(s_\lambda\) for the Schubert classes on \(G_3\) and
\(\tau_\lambda\) for those on \(G_6\). They are indexed by partitions
contained in a rectangle with three rows and seven columns, or six
rows and four columns, respectively.
For a partition \(\mu\) outside the latter rectangle, that is, with
\(\mu_1>4\) or \(\ell(\mu)>6\), we set \(\tau_\mu=0\).
Schubert classes indexed by one-row partitions are called \emph{special}
Schubert classes. The Giambelli formula expresses any Schubert class
as a determinant of special Schubert classes.
On \(G_6\), the special Schubert class \(\tau_r\) equals
\(c_r(\cQ)\). Using the same notation for the restriction of
\(\tau_\lambda\) to \(Y\), the formula becomes
\[
\tau_\lambda=\det\left(q_{\lambda_i-i+j}\right)_{1\le i,j\le\ell(\lambda)}.
\]
The entries \(q_j\) are the Chern classes introduced above, with the
same index conventions. In particular,
\(\tau_r=q_r\) and \(\tau_{1^r}=u_r\), where \(1^r\) denotes the
partition consisting of \(r\) ones. The class of the empty partition
is one. We omit trailing zeros when convenient, and also suppress
parentheses and commas; for example,
\(s_{750}=s_{(7,5,0)}=s_{(7,5)}\).

Let \(\cU_{G_3}\) be the rank-three universal subbundle on \(G_3\).
For a rank-\(r\) bundle \(B\), its pullback to the complete flag bundle
admits a filtration by subbundles with line bundle quotients.
The first Chern classes of these quotients are called the \emph{Chern roots}
of \(B\). Pullback to this flag bundle is injective on Chow groups,
so computations with symmetric polynomials in Chern roots can be
used to establish identities in the original Chow ring. This is the
splitting principle.

We pad partitions with zeros as needed. For \(\ell(\lambda)\le r\),
define the \emph{Schur polynomial} by
\[
S_\lambda(x_1,\ldots,x_r)
\colonequals \frac{\det\left(x_i^{\lambda_j+r-j}\right)_{i,j=1}^r}
{\det\left(x_i^{r-j}\right)_{i,j=1}^r}.
\]
Both determinants are alternating polynomials, and their quotient is
a symmetric polynomial. For a bundle \(B\), we write \(S_\lambda(B)\)
for the Chow class obtained by substituting its Chern roots into
\(S_\lambda\). The Giambelli formula on \(G_3\) gives
\(s_\lambda=S_\lambda(\cU_{G_3}^\vee)\).
The composite \(j\circ p\) sends \((K,W)\) to \(K\in G_3\), so
\(p^*j^*\cU_{G_3}\) is the rank-three universal subbundle on the
relative Grassmannian \(G_\sigma\). Naturality of Chern classes gives
\begin{equation}\label{eq:pullback}
p^*j^*s_\lambda=S_\lambda\left((p^*j^*\cU_{G_3})^\vee\right).
\end{equation}
From now on, we denote this bundle by
\(\cK\colonequals p^*j^*\cU_{G_3}\subset q^*\cU_Y\).

The notation \([x_1^{a_1}x_2^{a_2}x_3^{a_3}]\) means extraction of
the coefficient of the indicated monomial. We define the Chern
polynomial by \(c_t(B)\colonequals \sum_{j\ge0}c_j(B)t^j\), and the Segre
polynomial by \(s_t(B)\colonequals c_t(B)^{-1}\).

\begin{lemma}[Gysin formula for a Grassmann bundle, {\cite[Theorem~1.1, Proposition~4.5]{DarondeauPragacz}}]\label{lem:gysin}
Let \(x_1,x_2,x_3\) be the Chern roots of \(\cK^\vee\), and set
\(\Delta(x)=\prod_{i<j}\left(x_i-x_j\right)\). For a symmetric polynomial \(f\)
in \(x\), we have
\begin{equation}\label{eq:coefficient}
q_*f(x)=\left[x_1^5x_2^4x_3^3\right]
\left(f(x)\Delta(x)\prod_{i=1}^3\sum_{j=0}^4q_jx_i^{-j}\right).
\end{equation}
In particular, for any partition \(\lambda\) with \(\ell(\lambda)\le3\),
\begin{equation}\label{eq:schur-push}
q_*S_\lambda(\cK^\vee)=
\begin{cases}
\tau_{\lambda-(3,3,3)},&\lambda_3\ge3,\\
0,&\lambda_3<3.
\end{cases}
\end{equation}
\end{lemma}
\begin{proof}
Applying Theorem~1.1 of \cite{DarondeauPragacz} to \(\cU_Y\),
with \(n=6\) and \(d=3\), and using
\(s_t(\cU_Y)=c_t(\cQ_Y)\), gives \eqref{eq:coefficient}.
Proposition~4.5 gives
\[
q_*S_\lambda(\cK^\vee)
=\det\left(q_{\lambda_i-3-i+j}\right)_{i,j=1}^3.
\]
If \(\lambda_3<3\), the last row is zero. Otherwise, the Giambelli
formula gives \eqref{eq:schur-push}.
\end{proof}

\section{Chow relations on a Debarre--Voisin fourfold}
\label{sec:relations}

\subsection{A relation from the tangent bundle}

\begin{lemma}\label{lem:symplectic-relation}
The odd Chern classes of the tangent bundle \(T_Y\) vanish, and
\begin{equation}\label{eq:first}
3h^3-8hq_2+5q_3=0.
\end{equation}
\end{lemma}
\begin{proof}
The symplectic form on \(Y\) gives \(T_Y\simeq T_Y^\vee\).
Since \(c_j(T_Y^\vee)=(-1)^jc_j(T_Y)\), we have \(2c_j(T_Y)=0\)
for odd \(j\). We work with rational coefficients, so in particular
\(c_1(T_Y)=c_3(T_Y)=0\).

Let \(x_1,\ldots,x_6\) be the Chern roots of \(\cU_Y^\vee\).
We suppress pullbacks to the flag bundle in the notation.
Since \(u_j=c_j(\cU_Y^\vee)\) and \(u_1=h\), we have
\[
h=\sum_i x_i,\qquad
u_2=\sum_{i<j}x_ix_j,\qquad
u_3=\sum_{i<j<k}x_ix_jx_k.
\]

For a bundle \(B\) with Chern roots \(\xi_\alpha\), its Chern character
is \(\ch(B)=\sum_\alpha e^{\xi_\alpha}\), with formal exponential
series. It is additive in short exact sequences and multiplicative
under tensor products. The tautological sequence therefore gives
\[
\ch(\cQ_Y)=10-\ch(\cU_Y).
\]
The Chern roots of \(\cU_Y\) are \(-x_i\), and those of
\(\bigwedge^3\cU_Y^\vee\) are \(x_i+x_j+x_k\), for \(i<j<k\).
Thus the normal sequence \eqref{eq:normal} gives
\[
\ch(T_Y)=\ch\left(\cU_Y^\vee\right)\left(10-\ch(\cU_Y)\right)
-\ch\left(\bigwedge^3\cU_Y^\vee\right).
\]
Substituting the roots yields
\begin{equation}\label{eq:chern-character}
\ch(T_Y)=10\sum_{i=1}^6e^{x_i}
-\sum_{i,j=1}^6e^{x_i-x_j}
-\sum_{i<j<k}e^{x_i+x_j+x_k}.
\end{equation}
Write \(\ch_m\) for the component of codimension \(m\) of the Chern
character, and set \(p_m\colonequals\sum_i x_i^m\) for \(m\ge1\).
These power sums satisfy
\(\ch_m(\cU_Y^\vee)=p_m/m!\).
The Newton identities express them in terms of the Chern classes;
in the degrees needed here, they give
\begin{equation}\label{eq:newton-low}
p_1=h,\qquad p_2=h^2-2u_2,\qquad p_3=h^3-3hu_2+3u_3.
\end{equation}
Exchanging the two indices shows that
\(\sum_{i,j}\left(x_i-x_j\right)^3=0\).
Each \(x_i\) occurs in ten triples, and each pair \(\{i,j\}\) in
four triples. Counting these occurrences gives
\[
\sum_{i<j<k}\left(x_i+x_j+x_k\right)^3
=10p_3+12\left(hp_2-p_3\right)+6u_3.
\]
Thus the component of codimension three in \eqref{eq:chern-character} is
\[
\ch_3(T_Y)=2p_3-2hp_2-u_3=-2hu_2+5u_3.
\]
Since \(c_1(T_Y)=0\), we have \(c_3(T_Y)=2\ch_3(T_Y)\).
Using \eqref{eq:u-q}, we find
\begin{equation}\label{eq:c3}
c_3(T_Y)=10u_3-4hu_2=2\left(3h^3-8hq_2+5q_3\right).
\end{equation}
Together with \(c_3(T_Y)=0\), this proves \eqref{eq:first}.
\end{proof}

\subsection{Relations from the hyperplane section}
\label{sec:hyperplane}

For a smooth variety \(X\), define the subgroup of
\emph{homologically trivial cycles} by
\[
\CH^i_{\mathrm{hom}}(X)\colonequals
\Ker\left(\cl:\CH^i(X)\longrightarrow H^{2i}(X,\Q)\right).
\]
The following vanishing theorem holds for the Chow groups of \(F_\sigma\).

\begin{theorem}[Voisin, {\cite[Theorem 2.4]{VoisinGeneralized}}]\label{thm:vanishing}
Let \(F_\sigma\subset\Gr(3,10)\) be a smooth Pl\"ucker hyperplane
section. Then \(\CH_k(F_\sigma)_{\mathrm{hom}}=0\) for \(k<9\).
In particular,
\(\CH^{12}_{\mathrm{hom}}(F_\sigma)=\CH^{13}_{\mathrm{hom}}(F_\sigma)=0\).
\end{theorem}

\begin{lemma}\label{lem:hyperplane}
Let \(L=s_1\in\CH^1(G_3)\), and let \(m=12,13\).
If \(\kappa\in\CH^m(G_3)\) satisfies \(L\kappa=0\), then
\(j^*\kappa=0\) in \(\CH^m(F_\sigma)\).
\end{lemma}
\begin{proof}
We have \(j_*j^*\kappa=L\kappa=0\).
Under Poincar\'e duality, the dual of
\(j_*:H^{2m}(F_\sigma,\Q)\to H^{2m+2}(G_3,\Q)\) is
\[
j^*:H^{40-2m}(G_3,\Q)\longrightarrow H^{40-2m}(F_\sigma,\Q).
\]
Since \(40-2m=16,14<20=\dim F_\sigma\), this restriction is an
isomorphism by the weak Lefschetz theorem. It follows that \(j_*\)
is injective on cohomology, and hence \(j^*\kappa\) is homologically
trivial. Theorem~\ref{thm:vanishing} then implies that it vanishes
in the Chow group as well.
\end{proof}

\begin{lemma}\label{prop:two}
The following relations hold in \(\CH^*(Y)\):
\begin{align}
44hq_2-19h^3&=0,\label{eq:second}\\
h^2q_2-hq_3-2q_2^2+7q_4&=0.\label{eq:third}
\end{align}
\end{lemma}
\begin{proof}
Since \(q\) has relative dimension nine, \(q_*p^*\) lowers codimension
by nine. To obtain relations in codimensions three and four on \(Y\),
we therefore seek classes in the kernels of the multiplication maps
\[
L:\CH^{12}(G_3)\longrightarrow\CH^{13}(G_3),\qquad
L:\CH^{13}(G_3)\longrightarrow\CH^{14}(G_3).
\]
The \emph{Pieri rule} is \(Ls_\lambda=\sum_\nu s_\nu\), where \(\nu\)
ranges over the partitions obtained by adding one box to \(\lambda\)
within the rectangle with three rows and seven columns.
Define the following classes in codimensions twelve and thirteen:
\begin{align*}
\kappa_{12}={}&-3s_{750}+6s_{741}-5s_{732}+3s_{660}-3s_{651}-s_{642}\\
&+5s_{633}+4s_{552}-4s_{543}+4s_{444},\\
\kappa_{13}={}&2s_{751}-4s_{742}+5s_{733}-2s_{661}+2s_{652}
-s_{643}-s_{553}+s_{544}.
\end{align*}
The Pieri rule gives \(L\kappa_{12}=L\kappa_{13}=0\).
For example, the coefficient of \(s_{760}\) in
\(L\kappa_{12}\) is \(-3+3=0\); all other coefficients
vanish similarly.
Lemma~\ref{lem:hyperplane} therefore gives
\(j^*\kappa_{12}=j^*\kappa_{13}=0\).
Applying \(q_*p^*\) to \(j^*\kappa_{12}=j^*\kappa_{13}=0\) and using
\eqref{eq:pullback} and \eqref{eq:schur-push}, we obtain
\begin{align*}
5\tau_3-4\tau_{21}+4\tau_{111}&=0,\\
5\tau_4-\tau_{31}-\tau_{22}+\tau_{211}&=0.
\end{align*}
Substituting \(\tau_3=q_3\), \(\tau_{21}=hq_2-q_3\), and
\(\tau_{111}=h^3-2hq_2+q_3\) into the first relation yields
\(4h^3-12hq_2+13q_3=0\).
Subtracting thirteen times \eqref{eq:first} from five times this
relation gives \eqref{eq:second}.
Substituting
\[
\tau_{31}=hq_3-q_4,\quad \tau_{22}=q_2^2-hq_3,\quad
\tau_{211}=h^2q_2-hq_3-q_2^2+q_4
\]
into the second relation gives \eqref{eq:third}.
\end{proof}

\subsection{A nowhere-vanishing section and its push-forward}
\label{sec:section}

\begin{lemma}\label{lem:nowhere-zero}
The rank-twelve bundle
\[
\cE\colonequals \bigwedge^2\cK^\vee\otimes q^*\cQ_Y^\vee
\]
on \(G_\sigma=\Gr_Y(3,\cU_Y)\) admits a nowhere-vanishing section.
Consequently, \(c_{12}(\cE)=0\).
\end{lemma}
\begin{proof}
For a point \((K,W)\in G_\sigma\) and \(\bar v\in V/W\), choose
any representative \(v\in V\) of \(\bar v\). Define
\[
t_{K,W}(k_1\wedge k_2,\bar v)\colonequals \sigma(k_1,k_2,v),
\qquad k_1,k_2\in K,\quad\bar v\in V/W.
\]
Since \(K\subset W\) and \(\sigma|_W=0\), this value is independent
of the representative \(v\).
The bundle map \(\bigwedge^2\cK\otimes V\to\cO_{G_\sigma}\)
vanishes on \(\bigwedge^2\cK\otimes q^*\cU_Y\), so \(t\) defines
an algebraic section of \(\cE\).

If \(t_{K,W}=0\), then \(\sigma(K,K,V)=0\).
The normal sequence \eqref{eq:normal} shows that \(d\sigma_W\) is
surjective. Taking \(w_i\in K\) in \eqref{eq:differential}, all three
terms vanish. Hence the composite of the surjections \(d\sigma_W\)
and \(\bigwedge^3W^\vee\to\bigwedge^3K^\vee\) is zero.
Its target is nonzero, a contradiction. See also
\cite[Lemma 4.5]{DHOV} for this argument concerning the restriction
of the differential. Thus \(t\) is nowhere vanishing.
The cokernel of \(\cO_{G_\sigma}\xrightarrow{t}\cE\) is a
rank-eleven vector bundle. The Whitney formula therefore gives
\(c_{12}(\cE)=0\).
\end{proof}

Set \(a=c_1(\cK^\vee)\). On \(G_\sigma\), we abbreviate
\(q^*h\) and \(q^*q_j\) to \(h\) and \(q_j\), respectively.
Lemma~\ref{lem:nowhere-zero} gives
\begin{equation}\label{eq:theta-zero}
q_*\left((2a-h)c_{12}(\cE)\right)=0\quad\text{in }\CH^4(Y).
\end{equation}

\begin{lemma}\label{prop:top-push}
The following equality holds in \(\CH^4(Y)\):
\begin{equation}\label{eq:top-push}
q_*\left((2a-h)c_{12}(\cE)\right)
=10\left(h^4-5h^2q_2+6q_2^2+q_4\right).
\end{equation}
Consequently,
\begin{equation}\label{eq:fourth}
h^4-5h^2q_2+6q_2^2+q_4=0.
\end{equation}
\end{lemma}
\begin{proof}
Let \(x_1,x_2,x_3\) be the Chern roots of \(\cK^\vee\).
The Chern classes of \(\cQ_Y^\vee\) are \((-1)^jq_j\), so if we set
\[
P(v)\colonequals v^4-hv^3+q_2v^2-q_3v+q_4,\qquad
C(x)\colonequals \prod_{i<j}P\left(x_i+x_j\right),
\]
then \(C(x)=c_{12}(\cE)\).
Expand \(H(x)\colonequals \left(2x_1+2x_2+2x_3-h\right)C(x)\) in Schur polynomials:
\(H=\sum_\lambda b_\lambda S_\lambda\).
Multiplying both sides by the denominator \(\Delta(x)\), we can
extract each coefficient as
\begin{equation}\label{eq:schur-extraction}
b_\lambda=\left[x_1^{\lambda_1+2}x_2^{\lambda_2+1}x_3^{\lambda_3}\right]
\left(H(x)\Delta(x)\right).
\end{equation}
The strictly decreasing exponent sequences associated with distinct
partitions are different. Thus this extraction picks out only the
term corresponding to the identity permutation in the alternating
determinant for \(S_\lambda\).
We assign degrees \(1\), \(1\), and \(j\) to \(x_i\), \(h\), and
\(q_j\), respectively. Let \(\mu\) be a partition of length at
most three. Since \(H\) has total degree thirteen, only the eleven
terms with \(\lambda=\mu+(3,3,3)\) and \(|\mu|\le4\) survive
\eqref{eq:schur-push}. Their coefficients and push-forwards are listed
in Table~\ref{tab:schur}. The product expansions used to compute the
coefficients are given in Appendix~\ref{app:expansion}.

\begin{table}[htbp]
\centering
\caption{Schur coefficients contributing to the push-forward of \(H(x)\).}\label{tab:schur}
\renewcommand{\arraystretch}{1.18}
\begin{tabular}{@{}cll@{}}
\toprule
\(\lambda\)&\(b_\lambda\)&\(q_*S_\lambda\)\\
\midrule
733&\(70\)&\(q_4\)\\
643&\(170\)&\(hq_3-q_4\)\\
553&\(104\)&\(q_2^2-hq_3\)\\
544&\(80\)&\(h^2q_2-hq_3-q_2^2+q_4\)\\
\midrule
633&\(-190h\)&\(q_3\)\\
543&\(-296h\)&\(hq_2-q_3\)\\
444&\(-72h\)&\(h^3-2hq_2+q_3\)\\
\midrule
533&\(152h^2+108q_2\)&\(q_2\)\\
443&\(120h^2+84q_2\)&\(h^2-q_2\)\\
\midrule
433&\(-40h^3-124hq_2-56q_3\)&\(h\)\\
\midrule
333&\(\begin{gathered}2h^4+30h^2q_2+36hq_3\\[-2pt]+12q_2^2+30q_4\end{gathered}\)&\(1\)\\
\bottomrule
\end{tabular}
\end{table}

The third column is obtained from the determinant formula in the
proof of Lemma~\ref{lem:gysin}. For example, for \(\lambda=544\),
\[
q_*S_{544}=\det\begin{pmatrix}q_2&q_3&q_4\\1&h&q_2\\0&1&h\end{pmatrix}
=h^2q_2-hq_3-q_2^2+q_4.
\]
Adding the contributions within each of the five groups separated by
horizontal rules in the table gives, respectively,
\begin{align*}
&80h^2q_2-14hq_3+24q_2^2-20q_4,\\
&-72h^4-152h^2q_2+34hq_3,\\
&120h^4+116h^2q_2+24q_2^2,\\
&-40h^4-124h^2q_2-56hq_3,\\
&2h^4+30h^2q_2+36hq_3+12q_2^2+30q_4.
\end{align*}
Their sum is the right-hand side of \eqref{eq:top-push}.
Combining this with \eqref{eq:theta-zero} proves \eqref{eq:fourth}.
\end{proof}

\subsection{The ambient groups in codimensions three and four}

\begin{corollary}\label{cor:ambient-groups}\label{cor:degree4}
We have \(\CH^3_{\mathrm{amb}}(Y)=\Q h^3\) and
\(\CH^4_{\mathrm{amb}}(Y)=\Q h^4\), with
\begin{equation}\label{eq:degree3}
hq_2=\frac{19}{44}h^3,\qquad q_3=\frac1{11}h^3,
\end{equation}
and
\begin{equation}\label{eq:degree4}
h^2q_2=\frac{19}{44}h^4,\quad hq_3=\frac1{11}h^4,\quad
q_2^2=\frac{93}{484}h^4,\quad q_4=\frac3{484}h^4.
\end{equation}
\end{corollary}
\begin{proof}
Solving \eqref{eq:first} and \eqref{eq:second} gives \eqref{eq:degree3}.
By \eqref{eq:generators}, the ambient group in codimension three is
spanned by \(h^3\), which is nonzero because \(\int_Yh^4>0\).

Multiplying \eqref{eq:degree3} by \(h\) gives the first two equalities
in \eqref{eq:degree4}. Substituting them into \eqref{eq:third} and
\eqref{eq:fourth}, we obtain
\[
2q_2^2-7q_4=\frac{15}{44}h^4,\qquad
6q_2^2+q_4=\frac{51}{44}h^4.
\]
Solving this system gives the other two equalities.
By \eqref{eq:generators}, the ambient group in codimension four is
spanned by \(h^4\), which is nonzero because \(\int_Yh^4>0\).
\end{proof}

\section{The ambient Chow ring}
\label{sec:ambient-ring}

\begin{proof}[Proof of Theorem~\ref{thm:main}]
Set \(z=44q_2-19h^2\). By \eqref{eq:second} and \eqref{eq:degree4},
\[
hz=0,\qquad
z^2=1936q_2^2-1672h^2q_2+361h^4=11h^4.
\]
Hence there is a graded ring homomorphism
\(\Q[H,\zeta]/\left(H\zeta,\zeta^2-11H^4\right)\to\CH^*_{\mathrm{amb}}(Y)\).
It is surjective by \eqref{eq:generators} and
Corollary~\ref{cor:ambient-groups}.
In the source ring, \(H^5=H\zeta^2/11=0\), and
\(1,H,H^2,\zeta,H^3,H^4\) span the underlying vector space.

We check that the six corresponding classes are linearly independent in cohomology.
Classes in different codimensions are independent, and
\(1,h,h^3,h^4\) are nonzero.

It remains to prove that \(h^2,z\) are linearly independent in codimension two.
We use the same notation for their cohomology classes.
Since \(h\) is ample,
\[
d\colonequals \int_Y h^4>0.
\]
Suppose that
\[
ah^2+bz=0
\quad\text{in }H^4(Y,\Q),
\qquad a,b\in\Q.
\]
Multiplying this equality by \(h^2\) and \(z\), respectively, and integrating,
we obtain, using \(hz=0\) and \(z^2=11h^4\),
\[
\begin{pmatrix}0\\0\end{pmatrix}
=
\begin{pmatrix}
\displaystyle\int_Y h^2\left(ah^2+bz\right)\\[4pt]
\displaystyle\int_Y z\left(ah^2+bz\right)
\end{pmatrix}
=d\begin{pmatrix}1&0\\0&11\end{pmatrix}
\begin{pmatrix}a\\b\end{pmatrix}.
\]
Since \(d>0\), we have \(a=b=0\).
Thus \(h^2,z\) are linearly independent in cohomology.
It follows that the ring homomorphism is an isomorphism,
and the cycle class map is injective on the entire ambient ring.
\end{proof}

\section{Beauville--Voisin relations and the Franchetta property}
\label{sec:applications}

\subsection{Chern classes of the tangent bundle and Beauville--Voisin relations}

\begin{lemma}[Chern classes of the tangent bundle]\label{lem:tangent-chern}
The following equalities hold in \(\CH^*(Y)\):
\begin{align}
c_2(T_Y)&=8q_2-3h^2,\label{eq:c2}\\
c_4(T_Y)&=6h^4-34h^2q_2+hq_3+46q_2^2-4q_4.\label{eq:c4}
\end{align}
\end{lemma}
\begin{proof}
The equality \(c_2(T_Y)=5h^2-8u_2\), equivalent to \eqref{eq:c2},
already appears in \cite[proof of Lemma~4.5]{DebarreVoisin}.

Use the Chern roots \(x_1,\ldots,x_6\) and the power sums
\(p_m=\sum_i x_i^m\) introduced in the proof of
Lemma~\ref{lem:symplectic-relation}. In addition to
\eqref{eq:newton-low}, the Newton identities give
\[
p_4=h^4-4h^2u_2+4hu_3+2u_2^2-4u_4.
\]
The binomial expansion and the number of occurrences of each monomial
in the triples give
\begin{align*}
\sum_{i,j}\left(x_i-x_j\right)^2&=12p_2-2h^2,\\
\sum_{i<j<k}\left(x_i+x_j+x_k\right)^2&=6p_2+4h^2,\\
\sum_{i,j}\left(x_i-x_j\right)^4&=12p_4-8hp_3+6p_2^2,\\
\sum_{i<j<k}\left(x_i+x_j+x_k\right)^4&=6h^2p_2+6p_2^2+4hp_3-6p_4.
\end{align*}
Thus the components of codimensions two and four in
\eqref{eq:chern-character} are
\begin{align*}
\ch_2(T_Y)&=-5h^2+8u_2,\\
\ch_4(T_Y)&=-\tfrac5{12}h^4+\tfrac43h^2u_2
+\tfrac76hu_3-\tfrac53u_2^2-\tfrac23u_4.
\end{align*}
Since \(c_1(T_Y)=0\), we have
\(\ch_2(T_Y)=-c_2(T_Y)\) and
\(12\ch_4(T_Y)=c_2(T_Y)^2-2c_4(T_Y)\). Therefore
\begin{align*}
c_2(T_Y)&=5h^2-8u_2,\\
c_4(T_Y)&=15h^4-48h^2u_2-7hu_3+42u_2^2+4u_4.
\end{align*}
Substitution of \eqref{eq:u-q} gives \eqref{eq:c2} and \eqref{eq:c4}.
\end{proof}

\begin{lemma}[A consequence of Debarre--Voisin]\label{lem:picard}
A very general \(Y_\sigma\) has Picard number one.
\end{lemma}
\begin{proof}
For the inclusions \(j:F_\sigma\hookrightarrow G_3\) and
\(i:Y\hookrightarrow G_6\), set
\[
\begin{split}
K_F&\colonequals \Ker\left(j_*:H^{20}(F_\sigma,\Q)\to H^{22}(G_3,\Q)\right),\\
K_Y&\colonequals \Ker\left(i_*:H^2(Y,\Q)\to H^{42}(G_6,\Q)\right).
\end{split}
\]
A rational \((p,q)\)-class is a rational cohomology class whose
complexification lies in the \((p,q)\)-component of the Hodge decomposition.
By \cite[proof of Theorem~2.2(2), Corollary~2.7]{DebarreVoisin},
for very general \(\sigma\), the space \(K_F\) contains no nonzero
rational \((10,10)\)-class, and the incidence correspondence induces
an isomorphism \(q_*p^*:K_F\to K_Y\).
This isomorphism sends the \((r,s)\)-component to the
\((r-9,s-9)\)-component. Hence \(K_Y\) contains no nonzero rational
\((1,1)\)-class.

Corollary~2.7 also shows that \(i_*\) has rank one.
The projection formula gives
\(K_Y\subset\{\alpha\in H^2(Y,\Q)\mid\int_Y\alpha h^3=0\}\).
Both spaces have codimension one, so they are equal.
Since \(\int_Yh^4>0\), we have
\[
H^2(Y,\Q)=\Q h\oplus K_Y.
\]
Subtracting a suitable rational multiple of \(h\) from any rational
\((1,1)\)-class gives a rational \((1,1)\)-class in \(K_Y\), which
must be zero. Thus the rational \((1,1)\)-classes are exactly \(\Q h\),
and the Lefschetz \((1,1)\) theorem implies that the Picard number is one.
\end{proof}

\begin{proof}[Proof of Corollary~\ref{cor:bv}]
Using \eqref{eq:c2}, \eqref{eq:degree3}, and \eqref{eq:degree4}, we obtain
\[
hc_2(T_Y)=\frac5{11}h^3,\qquad
c_2(T_Y)^2=\frac{69}{121}h^4.
\]
Substituting \eqref{eq:degree4} into \eqref{eq:c4} gives
\[
c_4(T_Y)=
\left(6-34\frac{19}{44}+\frac1{11}
+46\frac{93}{484}-4\frac3{484}\right)h^4
=\frac{27}{121}h^4.
\]
These Chow relations do not require any assumption on the Picard number.

If the Picard number is one, then \(H^1(Y,\cO_Y)=0\) and \(h\ne0\)
imply \(\CH^1(Y)=\Pic(Y)\otimes\Q=\Q h\).
By Lemmas~\ref{lem:symplectic-relation} and \ref{lem:tangent-chern},
the Chern classes of the tangent bundle are ambient classes as well.
Hence
\(R^*(Y)=\Q[h,c_2(T_Y),c_4(T_Y)]\subset\CH^*_{\mathrm{amb}}(Y)\).
The restriction of the cycle class map to this subring is injective
by Theorem~\ref{thm:main}. By Lemma~\ref{lem:picard}, this conclusion
also applies to a very general fiber.
\end{proof}

\begin{remark}
Debarre and Voisin computed \(\int_Yh^4=1452\)
\cite[(11)]{DebarreVoisin}.
Setting \(o_Y\colonequals h^4/1452\in\CH^4(Y)\), we have \(\deg(o_Y)=1\),
and the relations established in the proof of Corollary~\ref{cor:bv}
take the form
\[
h^2c_2(T_Y)=660o_Y,\qquad c_2(T_Y)^2=828o_Y,
\qquad c_4(T_Y)=324o_Y.
\]
\end{remark}

\subsection{The Franchetta property for the universal family}

Let \(B^\circ\subset\mathbf P\left(\bigwedge^3V^\vee\right)\) be the nonempty
open subset parametrizing smooth four-dimensional zero schemes
\(Y_\sigma\), and let \(\mathcal Y\to B^\circ\) be the corresponding
universal family. For \(b\in B^\circ\), define the group of
\emph{generically defined cycles} by
\[
\mathrm{GD}\CH^i_{B^\circ}(Y_b)\colonequals 
\Imop\left(\CH^i(\mathcal Y)\to\CH^i(Y_b)\right).
\]
The family has the \emph{Franchetta property} in codimension \(i\)
if, for every \(b\in B^\circ\), the cycle class map is injective on
\(\mathrm{GD}\CH^i_{B^\circ}(Y_b)\) \cite{Franchetta}.

\begin{corollary}\label{cor:franchetta}
The above universal family has the Franchetta property in every
codimension. More precisely, for every \(b\in B^\circ\),
\[
\mathrm{GD}\CH^*_{B^\circ}(Y_b)=\CH^*_{\mathrm{amb}}(Y_b).
\]
\end{corollary}
\begin{proof}
Set \(A=\bigwedge^3V^\vee\). At each \(W\in G_6\), the evaluation
map sends \(\sigma\) to \(\sigma|_W\). An alternating three-form on
\(W\) extends to \(V\) after choosing a complement, so this map is
surjective. Thus its kernel \(\mathcal N\), defined by the exact sequence
\[
0\longrightarrow\mathcal N\longrightarrow A\otimes\cO_{G_6}
\longrightarrow\bigwedge^3\cU^\vee\longrightarrow0,
\]
is a vector bundle of rank
\(\binom{10}{3}-\binom63=120-20=100\).
The incidence variety over the full parameter space \(\mathbf P(A)\) is
\[
\overline{\mathcal Y}
=\{(W,[\sigma])\in G_6\times\mathbf P(A)\mid\sigma|_W=0\}
\simeq\mathbf P_{G_6}(\mathcal N).
\]
Here projectivization parametrizes lines. For fixed \(W\), the
condition is \(\C\sigma\subset\mathcal N_W\), so the fiber over
\(W\) is \(\mathbf P(\mathcal N_W)\simeq\mathbf P^{99}\).
Thus this projective bundle has relative dimension ninety-nine.

Write \(\rho:\overline{\mathcal Y}\to G_6\) and
\(\overline\pi:\overline{\mathcal Y}\to\mathbf P(A)\) for the
natural projections. The space
\(\mathcal Y=\overline\pi^{-1}(B^\circ)\) is an open subset of
\(\overline{\mathcal Y}\). The inclusion
\(\mathcal N\subset A\otimes\cO_{G_6}\) identifies the tautological
lines on both projectivizations with the same line \(\C\sigma\).
Consequently,
\[
\cO_{\mathbf P_{G_6}(\mathcal N)}(1)
\simeq\overline\pi^*\cO_{\mathbf P(A)}(1).
\]
Let \(\xi=c_1\left(\cO_{\mathbf P_{G_6}(\mathcal N)}(1)\right)\).
Since \(\overline\pi|_{Y_b}\) is constant, the restriction of this
line bundle to \(Y_b\) is trivial. Thus
\(\xi|_{Y_b}=0\) in \(\CH^1(Y_b)\).

Take any \(\gamma\in\CH^i(\mathcal Y)\).
By localization, the map
\(\CH^i(\overline{\mathcal Y})\to\CH^i(\mathcal Y)\) is surjective,
so \(\gamma\) has a lift
\(\overline\gamma\in\CH^i(\overline{\mathcal Y})\).
By the projective bundle formula, we can write
\[
\overline\gamma=\sum_{j=0}^{99}\rho^*\alpha_j\,\xi^j,
\qquad \alpha_j\in\CH^{i-j}(G_6),
\]
with Chow groups of negative codimension understood to be zero.
For the inclusion \(i_b:Y_b\hookrightarrow G_6\), we have
\(\rho|_{Y_b}=i_b\). Restricting to \(Y_b\) and using
\(\xi|_{Y_b}=0\) gives
\[
\gamma|_{Y_b}=i_b^*\alpha_0.
\]
Hence \(\mathrm{GD}\CH^i_{B^\circ}(Y_b)\subset
\CH^i_{\mathrm{amb}}(Y_b)\).
Conversely, the pullback of any \(\alpha\in\CH^i(G_6)\) along
\(\rho|_{\mathcal Y}\) restricts to \(i_b^*\alpha\) on \(Y_b\).
This proves the reverse inclusion.

By this equality, a generically defined class with zero cohomology
class is an ambient class with the same property. The injectivity in
Theorem~\ref{thm:main} implies that it is zero in \(\CH^i(Y_b)\).
The Franchetta property therefore holds in every codimension.
\end{proof}

\appendix
\section{Coefficients of the top Chern class}\label{app:expansion}

We give the expansions used to obtain the second column of
Table~\ref{tab:schur} in Lemma~\ref{prop:top-push}.
Assign weights \(1,2,3,4\) to the formal variables
\(h,q_2,q_3,q_4\), respectively, and weight zero to \(x_1,x_2,x_3\).
Let \(C_m\) be the component of \(C(x)\) of weight \(m\). Set
\[
a=x_1+x_2+x_3,\quad b=x_1x_2+x_1x_3+x_2x_3,\quad c=x_1x_2x_3,
\]
and
\[
\alpha=2a,\qquad \beta=a^2+b,\qquad\gamma=ab-c.
\]
The elementary symmetric polynomials in the three variables
\(x_1+x_2,x_1+x_3,x_2+x_3\) are \(\alpha,\beta,\gamma\).
Expanding \(C(x)=\prod_{i<j}P\left(x_i+x_j\right)\) through weight four gives
\begin{align*}
C_0={}&\gamma^4,\\
C_1={}&-h\beta\gamma^3,\\
C_2={}&q_2\beta^2\gamma^2+\left(h^2-2q_2\right)\alpha\gamma^3,\\
C_3={}&-q_3\beta^3\gamma+\left(3q_3-hq_2\right)\alpha\beta\gamma^2
+\left(-3q_3+3hq_2-h^3\right)\gamma^3,\\
C_4={}&q_4\left\{\beta^4-4\alpha\beta^2\gamma+\left(2\alpha^2+4\beta\right)\gamma^2\right\}\\
&+hq_3\left\{\alpha\beta^2\gamma-\left(2\alpha^2+\beta\right)\gamma^2\right\}
+q_2^2\left(\alpha^2-2\beta\right)\gamma^2+h^2q_2\beta\gamma^2.
\end{align*}
The weight-\(m\) component of \(H=(2a-h)C\) is
\(H_m=\alpha C_m-hC_{m-1}\), with \(C_{-1}=0\).
In the total-degree convention used in the main text, \(x_i,h,q_j\)
have degrees \(1,1,j\), respectively. Since \(H\) has total degree
thirteen, the degree of \(H_m\) in \(x\) is \(13-m\).
If \(\lambda_3\ge3\), then \(|\lambda|\ge9\). Thus, for \(m>4\),
the degree of \(H_m\) in \(x\) is less than nine, and no Schur
polynomial with \(\lambda_3\ge3\) occurs.
Thus the five expressions above, together with
\eqref{eq:schur-extraction}, give all entries in the table.

For example, only the weight-zero term \(H_0=2a\gamma^4\) contributes
to \(b_{733}\). Since
\(\gamma=(x_1+x_2)(x_1+x_3)(x_2+x_3)\),
\[
b_{733}=\left[x_1^9x_2^4x_3^3\right]\left(2a\gamma^4\Delta(x)\right)=70.
\]
The other ten coefficients are obtained in the same way.
Using these coefficients and the push-forward formula
\eqref{eq:schur-push}, summing the terms
\(b_\lambda q_*S_\lambda(\cK^\vee)\) gives the identity
\eqref{eq:top-push} in Lemma~\ref{prop:top-push}.

\end{document}